\documentclass[a4paper,12pt]{article}
\usepackage{amsmath, amsfonts, amssymb}
\usepackage[english, russian]{babel}
\usepackage {color}

\definecolor{mycolor1}{rgb}{0.020,0.388,0.757}
\definecolor{mycolor2}{rgb}{0.208,0.208,0.208}

\definecolor{mycolor1}{rgb}{0.020,0.388,0.757}
\definecolor{mycolor2}{rgb}{0.208,0.208,0.208}

\begin{document}
УДК 517.956

\begin{center}
Обобщенная задача Хольмгрена для эллиптического уравнения с несколькими
сингулярными коэффициентами
\end{center}

\begin{center}
The generalized Holmgren problem for elliptic equation with several singular coefficients
\end{center}

\begin{center}
Т.Г.Эргашев
\end{center}

\begin{center}
\textit{Институт математики им.В.И.Романовского АН РУз.}
\end{center}
\begin{center}
Узбекистан, Ташкент
\end{center}
\begin{center}
\textit{E-mail:}  \underline
{\textit{{ergashev.tukhtasin@gmail.com}}}
\end{center}

\textbf{Abstract.} Recently found all the fundamental solutions of a
multidimensional singular elliptic equation are expressed in terms of the
well-known Lauricella hypergeometric function in many variables. In this
paper, we find a unique solution of the generalized Holmgren problem for an elliptic equation with several
singular coefficients in explicit form. When finding a solution, we use
decomposition formulas and some adjacent relations for the Lauricella
hypergeometric function in many variables.



\bigskip

\section{Введение и постановка задач}

Рассмотрим уравнение
\begin{equation}
\label{eq1}
L_{\alpha} ^{(m,n)} \left( {u} \right): = {\sum\limits_{i = 1}^{m} {u_{x_{i}
x_{i}}} }   + {\sum\limits_{j = 1}^{n} {{\frac{{2\alpha _{j}}} {{x_{j}
}}}u_{x_{j}}} }   = 0
\end{equation}
в области $R_{m}^{n +}  = {\left\{ {x:x_{1} > 0,x_{2} > 0,...,x_{n} > 0}
\right\}},$ где $ x: = \left( {x_{1} ,...,x_{m}}  \right)$,\,$ m \ge 2,0 < n
\le m;$ \,$\alpha = \left( {\alpha _{1} ,...,\alpha _{n}}  \right),\alpha
_{j} $ -- действительные числа, причем $0 < 2\alpha _{j} < 1,j = \overline
{1,n} .$

Пусть $\Omega \subset R_{m}^{n +}  $ - первая $2^{ - n}$-ая часть
$m$-мерного шара радиуса $R$ c центром в начале координат:
\begin{equation*}
\label{eq2}
\Omega : = {\left\{ {x:\,x_{1}^{2} + ... + x_{m}^{2} < R^{2},\,\,x_{1} >
0,\,...,x_{n} > 0\,} \right\}}
\end{equation*}
и $S \subset R_{m}^{n +}  $ - такая же часть соответствующей этому шару
сферы:
\begin{equation*}
\label{eq3}
S: = {\left\{ {x:\,x_{1}^{2} + ... + x_{m}^{2} = R^{2},\,\,x_{1} >
0,\,...,x_{n} > 0\,} \right\}}.
\end{equation*}

Введем обозначения:
\begin{equation*}
\label{eq4}
\tilde {x}_{p} = \left( {x_{1} ,...,x_{p - 1} ,x_{p + 1} ,...,x_{m}}
\right) \in R_{m - 1} ;
\quad
{\rm {\bf 0}} = \left( {0,...,0} \right) \in R_{m - 1}
\end{equation*}
\begin{equation*}
\label{eq7}
x^{\left( {2\alpha}  \right)} = x_{1}^{2\alpha _{1}}\cdot  ...\cdot x_{n}^{2\alpha _{n}
} ;
\quad
\tilde {x}_{p}^{\left( {2\alpha}  \right)} = x_{1}^{2\alpha _{1}} \cdot ...\cdot x_{p -
1}^{2\alpha _{p - 1}}  x_{p + 1}^{2\alpha _{p + 1}} \cdot ...\cdot x_{n}^{2\alpha _{n}
} ;
\end{equation*}
\begin{equation*}
\label{eq8}
dx = dx_{1}\cdot ...\cdot dx_{m} ;
\quad
d\tilde {x}_{p} = dx_{1}\cdot ...\cdot dx_{p - 1} dx_{p + 1} \cdot ...\cdot dx_{m} ;
\end{equation*}
\begin{equation*}
\label{eq9}
D _{p} : = {\left\{ {x:\,x_{1}^{2} + ... + x_{p - 1}^{2} + x_{p +
1}^{2} + ... + x_{m}^{2} < R^{2},\,\,x_{1} > 0,\,...,x_{n} > 0\,}
\right\}};
\end{equation*}
\begin{equation*}
\label{eq10}
S_{p} : = {\left\{ {x:\,x_{1}^{2} + ... + x_{p - 1}^{2} + x_{p + 1}^{2} +
... + x_{m}^{2} = R^{2},\,\,x_{1} > 0,\,...,x_{n} > 0\,} \right\}}.
\end{equation*}

Во всех обозначениях $p = \overline {1,n} $.

\textbf{Задача Хольмгрена (Задача} $N$\textbf{).} Найти регулярное решение
$u_{0} \left( {x} \right)$ уравнения (\ref{eq1}) в области $\Omega $,
удовлетворяющее условиям
\begin{equation}
\label{eq11}
{\left. {\left( {x_{p}^{2\alpha _{p}}  {\frac{{\partial u_{0}}} {{\partial
x_{p}}} }} \right)} \right|}_{x_{p} = 0} = \nu _{p} \left( {\tilde {x}_{p}}
\right),\,\,\,\tilde {x}_{p} \in D _{p} ,\,\,p = \overline {1,n} ,
\quad
{\left. {u_{0}}  \right|}_{S} = \varphi _{0} \left( {x} \right),\,x \in S,\,
\end{equation}
где $\nu _{p} \in C\left( {D_{p}}  \right)$ и $\varphi _{0} \in C\left(
{\bar {S}} \right)$ - заданные функции, причем при стремлении точек $\tilde
{x}_{p} \in D _{p} $ к границе области $D _{p} $ функции $\nu _{p} $ могут
обращаться в бесконечность порядка меньше $1 - 2\alpha _{1} - ... - 2\alpha
_{n} $.

\textbf{Задача Дирихле.} Найти регулярное решение $u_{n} \left( {x} \right)$
уравнения (\ref{eq1}), удовлетворяющее условиям
\begin{equation}
\label{eq12}
{\left. {u_{n}}  \right|}_{x_{p} = 0} = \tau _{p} \left( {\tilde {x}_{p}}
\right),\,\,\tilde {x}_{p} \in \overline {D}  _{p} ,\,\,p = \overline
{1,n} ,
\quad
{\left. {u_{n}}  \right|}_{S} = \varphi _{n} \left( {x} \right),\,\,x \in
\overline {S} ,
\end{equation}
где $\tau _{p} \in C\left( {\bar {D} _{p}}  \right)$ и $\varphi _{n} \in
C\left( {\bar {S}} \right)$ - заданные функции, причем $\,\tau _{i} \left(
{{\rm {\bf 0}}} \right) = \tau _{j} \left( {{\rm {\bf 0}}} \right),\,i,j =
\overline {1,n} , \quad \,\tau _{i} \left( {\tilde {x}_{iq}^{0}}  \right) = \tau
_{j} \left( {\tilde {x}_{jq}^{0}}  \right),\,\,i,j,q = \overline {1,n}
,\,\,i \ne q,\,\,j \ne q, \quad {\left. {\tau _{p} \left( {\tilde {x}_{p}}
\right)} \right|}_{S_{p}}  = {\left. {\varphi _{n}}  \right|}_{S_{p}}
,\tilde {x}_{p} \in \overline {D}  _{p} $.

\textbf{Обобщенная задача Хольмгрена} (\textbf{Задача} $T^{k}N^{n -
k}$\textbf{).} Найти регулярное решение $u_{k} \left( {x} \right)$ уравнения
(\ref{eq1}), удовлетворяющее условиям
\begin{equation}
\label{eq13}
{\left. {u_{k}}  \right|}_{x_{p} = 0} = \tau _{p} \left( {\tilde {x}_{p}}
\right),\,\,\,p = \overline {1,k} ,
\end{equation}
\begin{equation}
\label{eq14}
{\left. {\left( {x_{p}^{2\alpha _{p}}  {\frac{{\partial u_{k}}} {{\partial
x_{p}}} }} \right)} \right|}_{x_{p} = 0} = \nu _{p} \left( {\tilde {x}_{p}}
\right),\,\,\,p = \overline {k + 1,n} ,
\end{equation}
\begin{equation}
\label{eq15}
{\left. {u_{k}}  \right|}_{S} = \varphi _{k} \left( {x} \right),
\end{equation}
где $\tau _{p} \in C\left( {\bar {D} _{p}}  \right)(p = \overline
{1,k} )$, $\nu _{p} \in C\left( {D _{p}}  \right)$ ($p = \overline {k +
1,n} )$ и $\phi _{k} \in C\left( {\bar {S}} \right)$ - заданные функции,
причем $\,\tau _{i} \left( {{\rm {\bf 0}}} \right) = \tau _{j} \left( {{\rm
{\bf 0}}} \right),\,i,j = \overline {1,n} , \quad \,\tau _{i} \left( {\tilde
{x}_{iq}^{0}}  \right) = \tau _{j} \left( {\tilde {x}_{jq}^{0}}  \right),$
$\,i \ne q,\,\,j \ne q, \quad {\left. {\tau _{p} \left( {\tilde {x}_{p}}
\right)} \right|}_{S_{p}}  = {\left. {\varphi _{k}}  \right|}_{S_{p}}
,\tilde {x}_{p} \in \overline {D}  _{p} , \quad i,j,p,q = \overline {1,k}
.$ Кроме того, при стремлении точек $\tilde {x}_{p} \in D _{p} $ к
границе области $ D _{p} $ функции $\nu _{p} $ могут обращаться в бесконечность
порядка меньше $1 - 2k + 2\alpha _{1} + ... + 2\alpha _{k} - 2\alpha _{k +
1} - ... - 2\alpha _{n} $.

Здесь и далее $k$ принимает неотрицательные целые значения: $k = 0,1,...,n.$

Фундаментальные решения уравнения (\ref{eq1}) в двумерном случае были известны и они успешно применены к решению основных краевых
задач для эллиптического уравнения с одной линией вырождения \cite{A1}, а для уравнения с двумя сингулярными коэффициентами  в работе \cite{A2} построены фундаментальные решения и изучены
задачи Хольмгрена, Дирихле, обобщенная задача Хольмгрена (задача
$H^{1}N^{1})$ в конечной \cite{A3} и в бесконечной \cite{{A4},{A5}} областях.

В трехмерном случае для уравнения (\ref{eq1}) с одной,
двумя и тремя сингулярными коэффициентами
построены фундаментальные решения \cite{{A6},{A7},{A8}} и найдены формулы решений различных
вариантов обобщенной задачи Хольмгрена \cite{{A9},{A10},{A11},{A12},{A13}}.

Когда размерность уравнения превышает три исследованию уравнения (\ref{eq1}) только
с одним или двумя сингулярными коэффициентами посвящены работы \cite{{A14},{A15},{A16},{A17},{A18}}.

Недавно в работе \cite{A19} построены все фундаментальные решения уравнения (\ref{eq1})
при $m \ge 2,0 < n \le m$ и с помощью одного из которых выписано явное решение  \cite{A20}
задачи Дирихле для уравнения (\ref{eq1}) с условиями (\ref{eq12}).

В настоящей работе единственное решение обобщенной задачи Хольмгрена для уравнения
(\ref{eq1}) с условиями (\ref{eq13})-(\ref{eq15}) находится в явном
виде.

\section{Предварительные сведения из теории гипергеометрической функции
Лауричелла многих переменных}

Гипергеометрическая функция Гаусса определяется формулой \cite[c.69]{A21}
\begin{equation*}
\label{eq16}
F(a,b;c;z) = F{\left( {{\begin{array}{*{20}c}
 {a,b;} \hfill \\
 {c;} \hfill \\
\end{array}} z} \right)} = {\sum\limits_{p = 0}^{\infty}  {{\frac{{\left(
{a} \right)_{p} \left( {b} \right)_{p}}} {{p!\left( {c} \right)_{p}
}}}z^{p}}} ,\,\,\,{\left| {z} \right|} < 1,\,
\end{equation*}
где $a,b$ и $c$ - комплексные постоянные, причем $\,\,c \ne 0, - 1, -
2,...$, а $\left( {\mu}  \right)_{\lambda}  $ означает символ Похгаммера:
\[
\left( {\mu}  \right)_{0} = 1,
\left( {\mu}  \right)_{\lambda}  = \mu \left( {\mu + 1} \right) \cdot ...
\cdot \left( {\mu + \lambda - 1} \right) = {\frac{{\Gamma \left( {\mu +
\lambda}  \right)}}{{\Gamma \left( {\mu}  \right)}}},
\]
$\Gamma \left( {\delta}  \right)$ -- известная гамма-функция Эйлера ($\delta \ne 0, - 1, - 2,....)$.

Гипергеометрическая функция Лауричелла многих переменных $F_{A}^{\left( {n} \right)}$ определяется
формулой \cite{A22}
\[
F_{A}^{\left( {n} \right)} \left( {a,b_1,...b_n;c_1,...,c_n;z_1,...,z_n} \right) = F_{A}^{(n)} {\left[
{{\begin{array}{*{20}c}
 {a,b_{1} ,...,b_{n};}  \hfill \\
 {c_{1} ,...,c_{n};}  \hfill \\
\end{array}} z_{1} ,...,z_{n}}  \right]}
\]
\begin{equation}
\label{eq19}
 = {\sum\limits_{p_{1} ,...,p_{n} = 0}^{\infty}  {{\displaystyle\frac{{\left( {a}
\right)_{p_{1} + ... + p_{n}}  \left( {b_{1}}  \right)_{p_{1}}  ...\left(
{b_{n}}  \right)_{p_{n}}} } {{\left( {c_{1}}  \right)_{p_{1}}  ...\left(
{c_{n}}  \right)_{p_{n}}} } }}} {\frac{{z_{1}^{p_{1}}} } {{p_{1}
!}}}...{\frac{{z_{n}^{p_{n}}} } {{p_{n} !}}},
\quad
{\sum\limits_{i = 1}^{n} {{\left| {z_{i}}  \right|} < 1,}}
\end{equation}
где
\[
c_{i} \ne 0, - 1, - 2,...,\,\,i = \overline {1,n} .
\]

Для исследования свойств гипергеометрической функции Лауричелла,
определенной формулой (\ref{eq19}), необходимы формулы разложения, позволяющие
представить гипергеометрическую функцию многих переменных через бесконечную
сумму произведений нескольких гипергеометрических функций с одной
переменной, а это, в свою очередь, облегчает процесс изучения свойств
функций многих переменных.  Например, гипергеометрическая функция Лауричелла
$F_{A}^{\left( {n} \right)} $, определенная формулой (\ref{eq19}), имеет разложение
\cite{A23}
\[
\begin{array}{l}
 F_{A}^{(n)} \left( a,b_{1} ,...,b_{n} ;c_{1} ,...,c_{n} ;z_{1}
,...,{z _{n}}  \right) \\
\\
 = {\sum\limits_{m_{2} ,...,m_{n} = 0}^{\infty}  {}} {\displaystyle\frac{{\left( {a}
\right)_{m_{2} + ... + m_{n}}  \left( {b_{1}}  \right)_{m_{2} + ... + m_{n}
} \left( {b_{2}}  \right)_{m_{2}}  ...\left( {b_{n}}  \right)_{m_{n}}
}}{{m_{2} !...m_{n} !\left( {c_{1}}  \right)_{m_{2} + ... + m_{n}}  \left(
{c_{2}}  \right)_{m_{2}}  ...\left( {c_{n}}  \right)_{m_{n}}} } }{z
}_{1}^{m_{2} + ... + m_{n}}  z_{2}^{m_{2}}  ...{z}_{n}^{m_{n}}  \\
 \end{array}
\]
\[
 \cdot F\left( a + m_{2} + ... + m_{n} ,b_{1} + m_{2} + ... + m_{n} ;c_{1}
+ m_{2} + ... + m_{n} ;{z_{1}}  \right)
\]
\begin{equation}
\label{eq23}
 \cdot F_{A}^{(n - 1)} \left[ {\begin{array}{*{20}c}
 {a + m_{2} + ... + m_{n} ,b_{2} + m_{2} ,...,b_{n} + m_{n} ;} \hfill \\
 {c_{2} + m_{2} ,...,c_{n} + m_{n} ;} \hfill \\
\end{array}} z _{2} ,...,{z_{n}}  \right],\,n=2,3,....
\end{equation}

Однако, из-за рекуррентности формулы (\ref{eq23}) могут возникать дополнительные
трудности в приложениях этого разложения. Дальнейшее изучение свойств гипергеометрической функции Лауричелла
$F_{A}^{\left( {n} \right)} $ показало, что формулу (\ref{eq23}) можно привести к более
удобному виду.

А именно справедлива

\textbf{Лемма 1}\cite{{A19},{A20}}. \textit{При $n = 1,2,3,...$ имеет место следующая формула
разложения:
\[
F_{A}^{(n)} \left( {a,b_{1} ,...,b_{n} ;c_{1} ,...,c_{n} ;z_1,...z_n}  \right)={\sum\limits_{{\mathop {m_{i,j} =
0}\limits_{(2 \le i \le j \le n)}}} ^{\infty}  {{\frac{{(a)_{A(n,n)}}}{{{\mathop {m_{2,2} ! \cdot \cdot \cdot m_{i,j} ! \cdot
\cdot \cdot m_{n,n} !}}}}}}}
\]
\begin{equation}
\label{eq24}
 \cdot {\prod\limits_{k = 1}^{n}
{{\left[ {{\frac{{(b_{k} )_{B(k,n)}}} {{(c_{k} )_{B(k,n)}}} }z
_{k}^{B(k,n)} F\left( {{\begin{array}{*{20}c}
 {a + A(k,n),b_{k} + B(k,n);} \hfill \\
 {c_{k} + B(k,n);} \hfill \\
\end{array}} z_{k}}  \right)} \right]}}} ,
\end{equation}
где
\begin{equation*}
\label{eq25}
A\left( {k,n} \right) = {\sum\limits_{i = 2}^{k + 1} {{\sum\limits_{j =
i}^{n} {m_{i,j}}} } } ,
B(k,n) = {\sum\limits_{i = 2}^{k} {m_{i,k} +}}  {\sum\limits_{i = k + 1}^{n}
{m_{k + 1,i}}}  .
\end{equation*}}

Лемма 1 доказывается методом математической индукции.

При исследовании краевых задач для эллиптического уравнения с сингулярными коэффициентами очень важна формула разложения, когда одна из переменных функции Лауричелла обращается в нуль.

\textbf{Следствие 1}.  \textit{При $n = 2,3,...$ и $l= 1,2,...,n$ имеет место следующая формула разложения:
\[
F_{A}^{(n - 1)} \left( {a,b_{1} ,...,b_{l - 1} ,b_{l + 1} ,...,b_{n} ;c_{1}
,...,c_{l-1},c_{l+1},...,c_{n} ;z_{1} ,...,z_{l - 1} ,z_{l + 1} ,...,z_{n}}  \right)
\]
\[
= {\sum\limits_{{\mathop {m_{i,j} =
0}\limits_{(2 \le i \le j \le n-1)}}} ^{\infty}  {{\frac{{(a)_{A(n-1,n-1)}}}{{{\mathop {m_{2,2} !\cdot \cdot \cdot m_{i,j} ! \cdot
\cdot \cdot m_{n-1,n-1} !}}} }}}}
{\prod\limits_{k = 1}^{l - 1} {{ {{\frac{{(b_{k} )_{B(k,n - 1)}
}}{{(c_{k} )_{B(k,n - 1)}}} } }{\prod\limits_{k = l + 1}^{n} {{ {{\frac{{(b_{k} )_{B(k - 1,n -
1)}}} {{(c_{k} )_{B(k - 1,n - 1)}}} }} }}} }}}  \\
\]
\[
{\cdot\prod\limits_{k = 1}^{l - 1}z_{k}^{B(k,n - 1)} {{ {F\left(
{{\begin{array}{*{20}c}
 {a + A(k,n - 1),b_{k} + B(k,n - 1);} \hfill \\
 {c_{k} + B(k,n - 1);} \hfill \\
\end{array}} z_{k}}  \right)} }}}  \\
\]
\begin{equation}
\label{eq27}
 \cdot {\prod\limits_{k = l + 1}^{n}z_{k}^{B(k - 1,n - 1)} {{ {F\left(
{{\begin{array}{*{20}c}
 {a + A(k - 1,n - 1),b_{k} + B(k - 1,n - 1);} \hfill \\
 {c_{k} + B(k - 1,n - 1);} \hfill \\
\end{array}} z_{k}}  \right)}}}} . \\
 \end{equation}}

\textbf{Лемма 2.} \textit{Пусть $a,b_{1} ,$\ldots , $b_{n} $ -- действительные
числа, причем $a \ne 0,\, - 1,\, - 2,...$и $a > b_{1} + ... + b_{n} .$ Тогда
при $n = 1,2,...$ справедлива следующая формула суммирования:
\[
{\sum\limits_{{\mathop {m_{i,j} =
0}\limits_{(2 \le i \le j \le n)}}} ^{\infty}  {{\frac{{(a)_{A(n,n)}}}{{{\mathop {m_{2,2} ! \cdot \cdot \cdot m_{i,j} ! \cdot
\cdot \cdot m_{n,n} !}}} }}}} {\prod\limits_{k = 1}^{n} {{\frac{{\left( {b_{k}}
\right)_{B(k,n)} \left( {a - b_{k}}  \right)_{A(k,n) - B(k,n)}}} {{\left(
{a} \right)_{A(k,n)}}} }}}
\]
\begin{equation}
\label{eq28}
 = \Gamma \left( {a - {\sum\limits_{k = 1}^{n}
{b_{k}}} }  \right){\frac{{\Gamma ^{n - 1}\left( {a}
\right)}}{{{\prod\limits_{k = 1}^{n} {\Gamma \left( {a - b_{k}}  \right)}
}}}}.
\end{equation}}

Лемма 2 доказывается методом математической индукции.

Нетрудно заметить, что формула (\ref{eq28}) является естественным обобщением известной формулы
суммирования для гипергеометрической функции Гаусса \cite[c.112]{A2}:
\begin{equation}
\label{eq17}
F(a,b;c;1) = {\frac{{\Gamma \left( {c} \right)\Gamma \left( {c - a - b}
\right)}}{{\Gamma \left( {c - a} \right)\Gamma \left( {c - b}
\right)}}},\,\,\,Re(c - a - b) > 0.
\end{equation}

\textbf{Лемма 3.}  \textit{Справедливо равенство
\[
\lim \limits_{\mathop {z_{k} \to 0,}\limits_{k = \overline {1,n}}}
z_{1}^{ - b_{1}}  ...z_{n}^{ - b_{n}}  F_{A}^{(n)} \left(
{a,b_{1} ,...,b_{n} ;c_{1} ,...,c_{n} ;1 - {\frac{{1}}{{z_{1}}} },...,1 -
{\frac{{1}}{{z_{n}}} }} \right)
\]
\begin{equation}
\label{eq30}
 = {\frac{{1}}{{\Gamma (a)}}}\Gamma \left( {a - {\sum\limits_{k = 1}^{n}
{b_{k}}} }  \right){\prod\limits_{k = 1}^{n} {{\frac{{\Gamma \left( {c_{k}}
\right)}}{{\Gamma \left( {c_{k} - b_{k}}  \right)}}}}} , \,\,\,\,n=1,2,....
\end{equation}}

Доказательство леммы 3 следует из лемм 1 и 2.

\section{Фундаментальные решения уравнения (\ref{eq1})}

При решении краевых задач важную роль играют фундаментальные решения
уравнения (\ref{eq1}). Известно \cite{A19}, что фундаментальные решения уравнения (\ref{eq1}) выражаются через гипергеометрическую функцию Лауричелла многих переменных:
\begin{equation}
\label{eq31}
q_{k} \left( {x,\xi}  \right) = \gamma _{k}  \prod\limits_{i=1}^k\left(x_i\xi_i\right)^{1-2\alpha_i} \cdot r^{ - 2\beta _{k}} F_{A}^{\left(
{n} \right)} {\left[ {{\begin{array}{*{20}c}
 {\beta _{k} ,1 - \alpha _{1} ,...,1 - \alpha _{k} ,\alpha _{k + 1}
,...,\alpha _{n} ;} \hfill \\
 {2 - 2\alpha _{1} ,...,2 - 2\alpha _{k} ,2\alpha _{k + 1} ,...,2\alpha _{n}
;} \hfill \\
\end{array}} \theta}  \right]},
\end{equation}
где
\[
\beta _{k} = {\frac{{m}}{{2}}} + k - 1 - {\sum\limits_{i = 1}^{k} {\alpha
_{i}}}   + {\sum\limits_{i = k + 1}^{n} {\alpha _{i}}}  ;
m \ge 2;
\]
\begin{equation}
\label{eq32}
\gamma _{k} = 2^{2\beta _{k} - m}{\frac{{\Gamma \left( {\beta _{k}}
\right)}}{{\pi ^{m / 2}}}}{\prod\limits_{i = k + 1}^{n} {{\frac{{\Gamma
\left( {\alpha _{i}}  \right)}}{{\Gamma \left( {2\alpha _{i}}  \right)}}}}
}{\prod\limits_{j = 1}^{k} {{\frac{{\Gamma \left( {1 - \alpha _{j}}
\right)}}{{\Gamma \left( {2 - 2\alpha _{j}}  \right)}}}}} ,
\quad
k = \overline {0,n} ;
\end{equation}
\begin{equation*}
\label{eq33}
\theta : = \left( {\theta _{1} ,...,\theta _{n}}  \right),
\theta _{i} = 1 - {\frac{{r_{i}^{2}}} {{r^{2}}}}, i=\overline{1,n}.
\end{equation*}

Здесь и в дальнейшем ${\sum\limits_{i = l + 1}^{l} {}} $ означает нуль, если $l = 0$
или $l = n$, и точно также ${\prod\limits_{i = l + 1}^{l} {}} $ равно
единице, если $l = 0$ или $l = n$.

Воспользовавшись формулой дифференцирования \cite{App}
\[
 {\frac{{\partial}} {{\partial z_{i}}} }F_{A}^{(n)} \left( {a,b_{1}
,...,b_{n} ;c_{1} ,...,c_{n} ;z_{1} ,...,z_{n}}  \right) = \\
\]
\begin{equation}
\label{eq34}
 = {\frac{{ab_{i}}} {{c_{i}}} }F_{A}^{(n)} {\left[
{{\begin{array}{*{20}c}
 {a + 1,b_{1} ,...,b_{i - 1} ,b_{i} + 1,b_{i + 1} ,...,b_{n} ;} \hfill \\
 {c_{1} ,...,c_{i - 1} ,c_{i} + 1,c_{i + 1} ,...,c_{n} ;} \hfill \\
\end{array}} z_{1} ,...,z_{n}}  \right]} \\
 \end{equation}
и смежным соотношением
\[
 {\sum\limits_{i = 1}^{n} {{\frac{{b_{i}}} {{c_{i}}} }F_{A}^{(n)} {\left[
{{\begin{array}{*{20}c}
 {a + 1,b_{1} ,...,b_{i - 1} ,b_{i} + 1,b_{i + 1} ,...,b_{n} ;} \hfill \\
 {c_{1} ,...,c_{i - 1} ,c_{i} + 1,c_{i + 1} ,...,c_{n} ;} \hfill \\
\end{array}} z_{1} ,...,z_{n}}  \right]}}}  \\
\]
\begin{equation}
\label{eq35}
 = F_{A}^{(n)} {\left[ {{\begin{array}{*{20}c}
 {a + 1,b_{1} ,...,b_{n} ;} \hfill \\
 {c_{1} ,...,c_{n} ;} \hfill \\
\end{array}} z_{1} ,...,z_{n}}  \right]} - F_{A}^{(n)} {\left[
{{\begin{array}{*{20}c}
 {a,b_{1} ,...,b_{n} ;} \hfill \\
 {c_{1} ,...,c_{n} ;} \hfill \\
\end{array}} z_{1} ,...,z_{n}}  \right]} \\
\end{equation}
вычислим производную по внешней нормали к границе области $\Omega$ формулой
\begin{equation}
\label{eq36}
{\frac{{\partial q_{k} \left( {x;\xi}  \right)}}{{\partial {\rm {\bf n}}}}}
= {\sum\limits_{i = 1}^{m} {{\frac{{\partial q_{k} \left( {x;\xi}
\right)}}{{\partial x_{i}}} }\cos \left( {x_{i} ,{\rm {\bf n}}} \right)}} .
\end{equation}

С этой целью для краткости изложения введем обозначения:
\begin{equation*}
\label{eq37}
F_{A}^{\left( {n} \right)} (\beta ;\theta ): = F_{A}^{\left( {n} \right)}
{\left[ {{\begin{array}{*{20}c}
 {\beta ,1 - \alpha _{1} ,...,1 - \alpha _{k} ,\alpha _{k + 1} ,...,\alpha
_{n} ;} \hfill \\
 {2 - 2\alpha _{1} ,...,2 - 2\alpha _{k} ,2\alpha _{k + 1} ,...,2\alpha _{n}
;} \hfill \\
\end{array}} \theta}  \right]},
\end{equation*}
\[
F_{(l \le k)}^{\left( {n} \right)} \left( {\theta}  \right):= F_{A}^{\left(
{n} \right)} {\left[ {{\begin{array}{*{20}c}
 {1 + \beta _{k} ,1 - \alpha _{1} ,...,1 - \alpha _{l - 1},} \hfill \\
 {2 - 2\alpha _{1} ,...,2 - 2\alpha _{l - 1},} \hfill
\\
\end{array}} }  \right.}
\]
\begin{equation*}
\label{eq38}
{\left. {{\begin{array}{*{20}c}
 {2 - \alpha _{l}
,1 - \alpha _{l + 1} ,...,1 - \alpha _{k} ,\alpha _{k + 1} ,...,\alpha _{n}
;} \hfill \\
 {3 - 2\alpha _{l} ,2 - 2\alpha
_{l + 1} ,...,2 - 2\alpha _{k} ,2\alpha _{k + 1} ,...,2\alpha _{n} ;} \hfill
\\
\end{array}} \theta}  \right]},
\end{equation*}
\[
F_{\left( {l > k} \right)}^{\left( {n} \right)} \left( {\theta}  \right): =
F_{A}^{\left( {n} \right)} {\left[ {{\begin{array}{*{20}c}
 {1 + \beta _{k} ,1 - \alpha _{1} ,...,1 - \alpha _{k},}
\hfill \\
 {2 - 2\alpha _{1} ,...,2 - 2\alpha _{k},} \hfill \\
\end{array}} }  \right.}
\]
\begin{equation*}
\label{eq39}
{\left. {{\begin{array}{*{20}c}
 {\alpha _{k + 1}
,...,\alpha _{l - 1} ,1 + \alpha _{l} ,\alpha _{l + 1} ,...,\alpha _{n} ;}
\hfill \\
 {2\alpha _{k + 1} ,...,2\alpha _{l
- 1} ,1 + 2\alpha _{l} ,2\alpha _{l + 1} ,...,2\alpha _{n} ;} \hfill \\
\end{array}} \theta}  \right]}
\end{equation*}
и подробно остановимся на вычисление ${{{\partial q_{k}}}/{{\partial x_{i}}}}$ при $1 \le i \le k.$ Действительно, используя
формулу дифференцирования (\ref{eq34}), получим
\[
{\frac{{\partial q_{k} \left( {x;\xi}  \right)}}{{\partial x_{i}}} } =
\gamma _{k} {\frac{{1 - 2\alpha _{i}}} {{x_{i}}} }\prod\limits_{i=1}^k\left(x_i\xi_i\right)^{1-2\alpha_i} \cdot r^{ - 2\beta _{k}} F_{A}^{\left(
{n} \right)} \left( {\beta _{k} ;\theta}  \right)
\]
\[
 - 2\beta _{k} \gamma _{k} \left( {x_{i} - \xi _{i}}\right)  \prod\limits_{i=1}^k\left(x_i\xi_i\right)^{1-2\alpha_i} \cdot r^{ - 2\beta _{k} -2}F_{A}^{\left( {n} \right)} \left( {\beta _{k} ;\theta}  \right)
\]
\[
- 2\beta_{k} \gamma _{k} \xi _{i} \prod\limits_{i=1}^k\left(x_i\xi_i\right)^{1-2\alpha_i} \cdot r^{ - 2\beta _{k} - 2}F_{\left( {i \le k} \right)}^{\left( {n}
\right)} \left( {\theta}  \right)
\]
\[
 - 2\beta _{k} \gamma _{k} \left( {x_{i} - \xi _{i}}  \right)\prod\limits_{i=1}^k\left(x_i\xi_i\right)^{1-2\alpha_i} \cdot r^{ - 2\beta _{k} -2}{\sum\limits_{s = 1}^{k} {{\frac{{1 - \alpha _{s}}} {{2 - 2\alpha _{s}
}}}}} F_{\left( {s \le k} \right)}^{\left( {n} \right)} \left( {\theta}
\right)
\]
\begin{equation}
\label{eq40}
 - 2\beta _{k} \gamma _{k} \left( {x_{i} - \xi _{i}}  \right)\prod\limits_{i=1}^k\left(x_i\xi_i\right)^{1-2\alpha_i} \cdot r^{ - 2\beta _{k} -2}{\sum\limits_{s = k + 1}^{n} {{\frac{{\alpha _{s}}} {{2\alpha _{s}}} }}
}F_{\left( {s > k} \right)}^{\left( {n} \right)} \left( {\theta}  \right).
\end{equation}

Отсюда в силу смежного соотношения (\ref{eq35}) будем иметь
\[
{\frac{{\partial q_{k}}} {{\partial x_{i}}} } = \gamma _{k} {\frac{{1 -
2\alpha _{i}}} {{x_{i}}} }\prod\limits_{i=1}^k\left(x_i\xi_i\right)^{1-2\alpha_i} \cdot r^{ - 2\beta _{k}} F_{A}^{\left( {n} \right)} \left( {\beta _{k};\theta}  \right)
\]
\[
- 2\beta _{k} \gamma _{k} \left( {x_{i} - \xi _{i}}  \right)\prod\limits_{i=1}^k\left(x_i\xi_i\right)^{1-2\alpha_i} \cdot r^{ - 2\beta _{k} -2}F_{A}^{\left( {n} \right)} \left( {1 + \beta _{k} ;\theta}  \right)
\]
\begin{equation}
\label{eq41}
-2\beta _{k} \gamma _{k} \xi _{i} \prod\limits_{i=1}^k\left(x_i\xi_i\right)^{1-2\alpha_i} \cdot r^{ - 2\beta _{k} - 2}F_{\left(
{i \le k} \right)}^{\left( {n} \right)} \left( {\theta}  \right),
1 \le i \le k.
\end{equation}

Аналогично вычисляются производные $\partial q_{k} / \partial x_{i} $ и по
другим переменным $x_{i} $:
\[
{\frac{{\partial q_{k}}} {{\partial x_{i}}} } = - 2\beta _{k} \gamma _{k}
\xi _{i} \prod\limits_{i=1}^k\left(x_i\xi_i\right)^{1-2\alpha_i} \cdot r^{ -
2\beta _{k} - 2}F_{\left( {i > k} \right)}^{\left( {n} \right)} \left(
{\theta}  \right)
\]
\begin{equation}
\label{eq42}
 - 2\beta _{k} \gamma _{k} \left( {x_{i} - \xi _{i}}  \right)\prod\limits_{i=1}^k\left(x_i\xi_i\right)^{1-2\alpha_i} \cdot r^{ - 2\beta _{k} -2}F_{A}^{\left( {n} \right)} \left( {1 + \beta _{k} ;\theta}  \right),
k + 1 \le i \le n;
\end{equation}
\begin{equation}
\label{eq43}
{\frac{{\partial q_{k}}} {{\partial x_{i}}} } = - 2\beta _{k} \gamma _{k}
\left( {x_{i} - \xi _{i}}  \right)\prod\limits_{i=1}^k\left(x_i\xi_i\right)^{1-2\alpha_i} \cdot r^{ - 2\beta _{k} - 2}F_{A}^{\left( {n} \right)} \left(
{1 + \beta _{k} ;\theta}  \right),
n + 1 \le i \le m.
\end{equation}

Подставляя (\ref{eq41})--(\ref{eq43}) в (\ref{eq36}), получим искомую производную
\[
{\frac{{\partial q_{k} \left( {x;\xi}  \right)}}{{\partial {\rm {\bf n}}}}}
= 2\beta _{k} \gamma _{k}\prod\limits_{i=1}^k\left(x_i\xi_i\right)^{1-2\alpha_i} \cdot r^{ - 2\beta _{k}}  F_{A}^{\left( {n} \right)}
\left( {1 + \beta _{k} ;\theta}  \right){\frac{{\partial}} {{\partial {\rm
{\bf n}}}}}{\left[ {\ln {\frac{{1}}{{r}}}} \right]} \\
\]
\[
 - 2\beta _{k} \gamma _{k} \prod\limits_{i=1}^k\left(x_i\xi_i\right)^{1-2\alpha_i} \cdot r^{ - 2\beta _{k} - 2}{\sum\limits_{i = 1}^{k} {\xi _{i} F_{\left({i \le k} \right)}^{\left( {n} \right)} \left( {\theta}  \right)\cos \left(
{{\rm {\bf n}};x_{i}}  \right)}}  \\
\]
\[
 - 2\beta _{k} \gamma _{k} \prod\limits_{i=1}^k\left(x_i\xi_i\right)^{1-2\alpha_i} \cdot r^{ - 2\beta _{k} - 2}{\sum\limits_{i = k + 1}^{n} {\xi _{i}
F_{\left( {i > k} \right)}^{\left( {n} \right)} \left( {\theta}
\right)\cos \left( {{\rm {\bf n}};x_{i}}  \right)}}  \\
\]
\begin{equation}
\label{eq44}
 + \gamma _{k} \prod\limits_{i=1}^k\left(x_i\xi_i\right)^{1-2\alpha_i} \cdot
r^{ - 2\beta _{k}} F_{A}^{\left( {n} \right)} \left( {\beta _{k} ;\theta}
\right){\sum\limits_{i = 1}^{k} {{\frac{{1 - 2\alpha _{i}}} {{x_{i}}} }\cos
\left( {{\rm {\bf n}};x_{i}}  \right)}} . \\
\end{equation}

Имеет место

\textbf{Лемма 4}. \textit{Фундаментальные решения
уравнения (\ref{eq1}) обладают следующими свойствами:
\begin{equation*}
\label{eq45}
{\left. {q_{k} \left( {x,\xi}  \right)} \right|}_{x_{p} = 0} = 0,\,\,p =
\overline {1,k} ;
\quad
{\left. {\left( {x_{p}^{2\alpha _{p}}  {\frac{{\partial q_{k}}} {{\partial
x_{p}}} }} \right)} \right|}_{x_{p} = 0} = 0,\,\,\,p = \overline {k + 1,\,n}
.
\end{equation*}}

Доказательство Леммы 4 непосредственно следует из формул (\ref{eq31}) и (\ref{eq44}).

\section{Формула Грина и теорема единственности.}

Рассмотрим тождество
\begin{equation*}
\label{eq46}
x^{\left( {2\alpha}  \right)}{\left[ {uL_{\alpha} ^{(m,n)} \left( {w}
\right) - wL_{\alpha} ^{(m,n)} \left( {u} \right)} \right]} =
{\sum\limits_{i = 1}^{m} {{\frac{{\partial}} {{\partial x_{i}}} }{\left[
{x^{\left( {2\alpha}  \right)}\left( {u{\frac{{\partial w}}{{\partial x_{i}
}}} - w{\frac{{\partial u}}{{\partial x_{i}}} }} \right)} \right]}}} ,\,\,n
\le m.
\end{equation*}

Интегрируя обе части этого тождества по области $\Omega $ и используя
известную формулу Гаусса-Остроградского, получаем формулу Грина для уравнения (\ref{eq1})
\[
{\int_{\Omega}  {x^{\left( {2\alpha}  \right)}{\left[ {uL_{\alpha} ^{(m,n)}
\left( {w} \right) - wL_{\alpha} ^{(m,n)} \left( {u} \right)} \right]}dx}}
\]
\begin{equation}
\label{eq47}
= {\int_{D}  {x^{\left( {2\alpha}  \right)}{\sum\limits_{i = 1}^{m}
{{\left[ {\left( {u{\frac{{\partial w}}{{\partial x_{p}}} } -
w{\frac{{\partial u}}{{\partial x_{i}}} }} \right)\cos \left( {{\rm {\bf
n}},x_{i}}  \right)} \right]}dD}} } } ,
\end{equation}
где $D$  -- граница области $\Omega$.

\textbf{Теорема 1}. \textit{Обобщенная задача Хольмгрена для уравнения (\ref{eq1}) в области $\Omega $  имеет не более одного решения. }

\textit{Доказательство.}
Поскольку $D $ -- граница области $\Omega $ и ${\rm {\bf n}}$ -- внешняя
нормаль к $D $, то нетрудно убедиться в том, что  $\cos \left( {{\rm {\bf n}},x_{i}}  \right)dD_{j}=d\tilde {x}_{i},$ если\,\,$i = j$ и $\cos \left( {{\rm {\bf n}},x_{i}}  \right)dD_{j}=0 $ в противном случае, где  $i = \overline {1,m} ,\,j = \overline {1,n}$.

Легко проверить равенство:
\begin{equation}
\label{eq49}
{\int_{\Omega}  {x^{\left( {2\alpha}  \right)}uL_{\alpha} ^{(m,n)}
\left( {u} \right)dx}}  = - {\int_{\Omega}  {x^{\left( {2\alpha}
\right)}{\sum\limits_{i = 1}^{m} {\left( {{\frac{{\partial u}}{{\partial
x_{i}}} }} \right)^{2}}} dx}}  + {\sum\limits_{i = 1}^{m}
{{\int_{\Omega}  {{\frac{{\partial}} {{\partial x_{i}}} }\left(
{x^{\left( {2\alpha}  \right)}u{\frac{{\partial u}}{{\partial x_{i}}} }}
\right)dx}}} } .
\end{equation}

После применения формулы Гаусса-Остроградского к интегралам, входящим в (\ref{eq49}), при условии, что $u(x) =
u_{k} (x)$ является решением уравнения (\ref{eq1}), получим
\[
{\int_{\Omega}  {x^{\left( {2\alpha}  \right)}{\sum\limits_{i = 1}^{m}
{\left( {{\frac{{\partial u_{k}}} {{\partial x_{i}}} }} \right)^{2}}} dx}}   = {\int_{S} {x^{\left( {2\alpha}  \right)}\varphi _{k} \left( {S}
\right){\frac{{\partial u_{k}}} {{\partial {\rm {\bf n}}}}}dS}}
\]
\begin{equation}
\label{eq50}
 + {\sum\limits_{i = 1}^{k} {{\int_{D _{i}}  {\tilde {x}_{i}^{\left(
{2\alpha}  \right)} \tau _{i} \left( {\tilde {x}_{i}}  \right)f_{i} \left(
{\tilde {x}_{i}}  \right)d\tilde {x}_{i}}} } }  + {\sum\limits_{i = k +
1}^{n} {{\int_{D _{i}}  {\tilde {x}_{i}^{\left( {2\alpha}  \right)}
\tilde {\nu} _{i} \left( {\tilde {x}_{i}}  \right)h_{i} \left( {\tilde
{x}_{i}}  \right)d\tilde {x}_{i}}} } } ,
\end{equation}
где
\[
f_{i} \left( {\tilde {x}_{i}}  \right): = {\mathop {\lim} \limits_{x_{i} \to
0}} x_{i}^{2\alpha _{i}}  {\frac{{\partial u_{k} (x)}}{{\partial x_{i}
}}},\,\,\,\,i = \overline {1,k} ,
\quad
h_{i} \left( {\tilde {x}_{i}}  \right): = {\mathop {\lim} \limits_{x_{i} \to
0}} u_{k} (x),\,\,\,\,i = \overline {k + 1,\,n} .
\]

Если рассмотреть однородный случай обобщенной задачи Хольмгрена, то из (\ref{eq50})
получается
\begin{equation*}
\label{eq51}
{\int_{\Omega}  {x^{\left( {2\alpha}  \right)}{\sum\limits_{i = 1}^{m}
{\left( {{\frac{{\partial u_{k}}} {{\partial x_{i}}} }} \right)^{2}}} dx}}   = 0.
\end{equation*}

Отсюда следует, что $u_{k} \left( {x} \right) = 0$ в $\overline {\Omega}  .$
Тем самым доказана единственность решения обобщенной задачи Хольмгрена для
уравнения (\ref{eq1}). Теорема 1 доказана.

\section{Решение обобщенной задачи Хольмгрена с помощью функции Грина.}

Существование решения докажем методом функции Грина. Для определенности положим $m>2$.

\textbf{Определение.} \textit{Функцией Грина обобщенной задачи Хольмгрена для
уравнения (\ref{eq1}) называется функция $G_{k} \left( {x;\xi}  \right)$,
удовлетворяющая условиям:}

(i)\,\textit{внутри области $\Omega $, кроме точки $\xi $, эта функция является
регулярным решением уравнения (\ref{eq1});}

(ii)\,\textit{удовлетворяет граничным условиям}

\begin{equation*}
\label{eq52}
{\left. {G_{k}}  \right|}_{x_{p} = 0} = 0,\,\,\,p = \overline {1,k} ,
\quad
{\left. {\left( {x_{p}^{2\alpha _{p}}  {\frac{{\partial G_{k}}} {{\partial
x_{p}}} }} \right)} \right|}_{x_{p} = 0} = 0,\,\,\,p = \overline {k + 1,n} ,
\quad
{\left. {G_{k}}  \right|}_{S} = 0,
\end{equation*}

(iii)\,\textit{она может быть представлена в виде
\begin{equation*}
\label{eq53}
G_{k} \left( {x;\xi}  \right) = q_{k} \left( {x;\xi}  \right) + g_{k} \left(
{x;\xi}  \right),
\end{equation*}
где $q_{k} \left( {x;\xi}  \right)$ -- фундаментальное решение уравнения (\ref{eq1}),
определенное формулой (\ref{eq31}), а функция $g_{k} \left( {x;\xi}  \right)$
является регулярным решением уравнения (\ref{eq1}) в области $\Omega $.}

Из представления функции Грина нам необходимо найти ее регулярную часть
$g_{k} \left( {x;\xi}  \right)$, удовлетворяющую условиям
\[
{\left. {g_{k} \left( {x;\xi}  \right)} \right|}_{S} = {\left. { - q_{k}
\left( {x;\xi}  \right)} \right|}_{S} ,
\quad
{\left. {g_{k} \left( {x;\xi}  \right)} \right|}_{x_{p} = 0} = - {\left.
{q_{k} \left( {x;\xi}  \right)} \right|}_{x_{p} = 0} ,\,\,\,p = \overline
{1,k} ,
\]

\[
{\left. {\left( {x_{p}^{2\alpha _{p}}  {\frac{{\partial g_{k} \left( {x;\xi
} \right)}}{{\partial x_{p}}} }} \right)} \right|}_{x_{p} = 0} = - {\left.
{\left( {x_{p}^{2\alpha _{p}}  {\frac{{\partial q_{k} \left( {x;\xi}
\right)}}{{\partial x_{p}}} }} \right)} \right|}_{x_{p} = 0} ,\,\,\,p =
\overline {k + 1,n} .
\]

Отсюда следует, что нужно выбрать регулярное решение $g_{k} \left( {x;\xi}  \right)$ в виде
\begin{equation*}
\label{eq541}
g_{k} \left( {x;\xi}  \right) =  - \left(
{{\frac{{R}}{{\rho}} }} \right)^{2\beta _{k}} q_{k} \left( {x;\bar {\xi}}
\right),
\end{equation*}
где
\[
\rho ^{2} = {\sum\limits_{p = 1}^{m} {\xi _{p}^{2}}}  ;
\quad
\bar {\xi}  = \left( {\bar {\xi} _{1} ,...,\bar {\xi} _{m}}  \right),
\bar {\xi} _{p} = {\frac{{R^{2}}}{{\rho ^{2}}}}\xi _{p} ,\,\,\,p = \overline
{1,m} .
\]

Следовательно, для области $\Omega$ функция Грина обобщенной задачи Хольмгрена для уравнения (\ref{eq1}) имеет вид:
\begin{equation*}
\label{eq54}
G_{k} \left( {x;\xi}  \right) = q_{k} \left( {x;\xi}  \right) - \left(
{{\frac{{R}}{{\rho}} }} \right)^{2\beta _{k}} q_{k} \left( {x;\bar {\xi}}
\right).
\end{equation*}

В дальнейших исследованиях нам будет необходимо вычислить производные по внешней
нормали к границе области $\Omega$ от функции $q_{k} \left( {x;\bar {\xi}}  \right)$:
\[
 {\frac{{\partial q_{k} \left( {x;\bar {\xi}}  \right)}}{{\partial {\rm {\bf
n}}}}} = 2\beta _{k} \gamma _{k}\prod\limits_{i=1}^k\left(x_i \bar{\xi_i}\right)^{1-2\alpha_i} \cdot \bar {r}^{ - 2\beta _{k}} F_{A}^{\left( {n}\right)} \left( {1 + \beta _{k} ;\bar {\theta}}  \right){\frac{{\partial
}}{{\partial {\rm {\bf n}}}}}{\left[ {\ln {\frac{{1}}{{\bar {r}}}}} \right]}
\\
\]
\[
 - 2\beta _{k} \gamma _{k} \prod\limits_{i=1}^k\left(x_i \bar{\xi_i}\right)^{1-2\alpha_i} \cdot \bar {r}^{ - 2\beta _{k} - 2}{\sum\limits_{i = 1}^{k}{\bar {\xi} _{i} F_{\left( {i \le k} \right)}^{\left( {n} \right)} \left(
{\bar {\theta}}  \right)\cos \left( {{\rm {\bf n}};x_{i}}  \right)}}  \\
\]
\[
 - 2\beta _{k} \gamma _{k} \prod\limits_{i=1}^k\left(x_i \bar{\xi_i}\right)^{1-2\alpha_i} \cdot \bar {r}^{ - 2\beta _{k} - 2}{\sum\limits_{i = k + 1}^{n}{\bar {\xi} _{i} F_{\left( {i > k} \right)}^{\left( {n} \right)} \left(
{\bar {\theta}}  \right)\cos \left( {{\rm {\bf n}};x_{i}}  \right)}}  \\
\]
\begin{equation}
\label{eq55}
 + \gamma _{k} \prod\limits_{i=1}^k\left(x_i \bar{\xi_i}\right)^{1-2\alpha_i} \cdot \bar {r}^{ - 2\beta _{k}} F_{A}^{\left( {n} \right)} \left( {\beta_{k} ;\bar {\theta}}  \right){\sum\limits_{i = 1}^{k} {{\frac{{1 - 2\alpha
_{i}}} {{x_{i}}} }\cos \left( {{\rm {\bf n}};x_{i}}  \right)}} , \\
\end{equation}
где
\begin{equation*}
\label{eq56}
\bar {\theta} : = \left( {\bar {\theta} _{1} ,...,\bar {\theta} _{n}}
\right),
\,\,\bar {\theta} _{i} = 1 - {\frac{{\bar {r}_{i}^{2}}} {{\bar {r}^{2}}}};   \,\,\bar {r}^{2} = {\sum\limits_{j = 1}^{m} {\left( {x_{j} -
{\frac{{R^{2}}}{{\rho ^{2}}}}\xi _{j}}  \right)^{2}}} ,
\end{equation*}
\begin{equation*}
\label{eq57}
\bar {r}_{i}^{2} = \left( {x_{i} + {\frac{{R^{2}}}{{\rho ^{2}}}}\xi _{i}}
\right)^{2} + {\sum\limits_{j = 1,j \ne i}^{m} {\left( {x_{j} -
{\frac{{R^{2}}}{{\rho ^{2}}}}\xi _{j}}  \right)^{2}}} ,
i = \overline {1,n} .
\end{equation*}

Пусть $\xi \in \Omega $. Вырежем из области $\Omega $ шар малого радиуса
$\varepsilon $ с центром в точке $\xi $ и оставшуюся часть $\Omega $
обозначим через $\Omega _{\varepsilon}  $, а через $C_{\varepsilon}  \, - $
сферу вырезанного шара. Используя формулу (\ref{eq47}), получим
\[
{\int_{C_{\varepsilon}}   {x^{\left( {2\alpha}  \right)}{\left[ {u_{k}
{\frac{{\partial G_{k}}} {{\partial {\rm {\bf n}}}}} - G_{k}
{\frac{{\partial u_{k}}} {{\partial {\rm {\bf n}}}}}}
\right]}dC_{\varepsilon}} }   = {\sum\limits_{p = 1}^{k} {{\int_{D _{p}
} {\tau _{p} \left( {\tilde {x}_{p}}  \right)\tilde {G}_{k} \left(
{x_{p}^{0} ;\xi}  \right)d\tilde {x}_{p}}} } }  \\
\]
\begin{equation}
\label{eq58}
 - {\sum\limits_{p = k + 1}^{n} {{\int_{D _{p}}  {\nu _{p} \left(
{\tilde {x}_{p}}  \right)\tilde {x}_{p}^{\left( {2\alpha}  \right)} G_{k}
\left( {x_{p}^{0} ;\xi}  \right)d\tilde {x}_{p}}} } }  - {\int_{S} {\varphi
_{k} \left( {S} \right)x^{\left( {2\alpha}  \right)}{\frac{{\partial G_{k}
\left( {x;\xi}  \right)}}{{\partial {\rm {\bf n}}}}}dS}} , \\
\end{equation}
где
\begin{equation*}
\label{eq59}
\tilde {G}_{k} \left( {x_{p}^{0} ;\xi}  \right) = \tilde {x}_{p}^{\left(
{2\alpha}  \right)} {\left. {\left( {x_{p}^{2\alpha _{p}}  {\frac{{\partial
G_{k} \left( {x,\xi}  \right)}}{{\partial x_{p}}} }} \right)}
\right|}_{x_{p} = 0} ,\,\,\,\,p= \overline {1,k};
\end{equation*}
\begin{equation*}
\label{eq60}
x_{p}^{0} = \left( {x_{1} ,...,x_{p - 1} ,0,x_{p + 1} ,...,x_{m}}  \right)
\in R_{m} ,\,\,p = \overline {1,n} .
\end{equation*}

В равенстве (\ref{eq58}) совершим предельный переход при $\varepsilon \to 0.$
Предварительно преобразуем левую часть (\ref{eq58})
\begin{equation*}
\label{eq61}
I = {\int_{C_{\varepsilon}}   {x^{\left( {2\alpha}  \right)}{\left[ {u_{k}
{\frac{{\partial G_{k}}} {{\partial {\rm {\bf n}}}}} - G_{k}
{\frac{{\partial u_{k}}} {{\partial {\rm {\bf n}}}}}}
\right]}dC_{\varepsilon}} }   = I_{1} - I_{2} ,
\end{equation*}
где
\[
I_{1} = {\int_{C_{\varepsilon}}   {x^{\left( {2\alpha}  \right)}u_{k} \left(
{x} \right){\frac{{\partial G_{k} \left( {x;\xi}  \right)}}{{\partial {\rm
{\bf n}}}}}dC_{\varepsilon}} }  ,
\quad
I_{2} = {\int_{C_{\varepsilon}}   {x^{\left( {2\alpha}  \right)}G_{k} \left(
{x;\xi}  \right){\frac{{\partial u_{k} \left( {x} \right)}}{{\partial {\rm
{\bf n}}}}}dC_{\varepsilon}} }  .
\]

Рассмотрим интеграл
\begin{equation*}
\label{eq62}
I_{1} = {\int_{C_{\varepsilon}}   {x^{\left( {2\alpha}  \right)}u_{k} \left(
{x} \right){\left[ {{\frac{{\partial q_{k} \left( {x;\xi}
\right)}}{{\partial {\rm {\bf n}}}}} - \left( {{\frac{{R}}{{\rho}} }}
\right)^{2\beta _{k}} {\frac{{\partial q_{k} \left( {x;\bar{\xi}}
\right)}}{{\partial {\rm {\bf n}}}}}} \right]}dC_{\varepsilon}} }   = I_{3}
- I_{4} .
\end{equation*}

Перепишем интеграл $I_{3} $ в виде
\begin{equation*}
\label{eq63}
I_{3} = {\int_{C_{\varepsilon}}   {u_{k} \left( {x} \right)x^{\left(
{2\alpha}  \right)}{\frac{{\partial q_{k} \left( {x;\xi}
\right)}}{{\partial {\rm {\bf n}}}}}dC_{\varepsilon}} }   = I_{5} +
I_{6} + I_{7} + I_{8} ,
\end{equation*}
где
\begin{equation*}
\label{eq64}
I_{5} = 2\beta _{k} \gamma _{k} {\int_{C_{\varepsilon}}   {u_{k} \left(
{x} \right)x^{\left( {2\alpha}  \right)}\prod\limits_{i=1}^k\left(x_i\xi_i\right)^{1-2\alpha_i} \cdot r^{ - 2\beta _{k}}}
}F_{A}^{\left( {n} \right)} \left( {1 + \beta _{k} ;\theta}
\right){\frac{{\partial}} {{\partial {\rm {\bf n}}}}}{\left[ {\ln
{\frac{{1}}{{r}}}} \right]}dC_{\varepsilon}  ,
\end{equation*}
\[
I_{6} = - 2\beta _{k} \gamma _{k} {\int_{C_{\varepsilon}}   {u_{k}
(x)x^{\left( {2\alpha}  \right)}\prod\limits_{i=1}^k\left(x_i\xi_i\right)^{1-2\alpha_i} \cdot r^{ - 2\beta _{k} - 2}{\sum\limits_{i = 1}^{k} {\xi_{i} F_{\left( {i \le k} \right)}^{\left( {n} \right)} \left( {\theta}
\right)\cos \left( {{\rm {\bf n}};x_{i}}  \right)}} dC_{\varepsilon}} }  ,
\]
\[
I_{7} = - 2\beta _{k} \gamma _{k} {\int_{C_{\varepsilon}}   {u_{k} \left(
{x} \right)x^{\left( {2\alpha}  \right)}\prod\limits_{i=1}^k\left(x_i\xi_i\right)^{1-2\alpha_i} \cdot r^{ - 2\beta _{k} -
2}{\sum\limits_{i = k + 1}^{n} {\xi _{i} F_{\left( {i > k} \right)}^{\left(
{n} \right)} \left( {\theta}  \right)\cos \left( {{\rm {\bf n}};x_{i}}
\right)}} dC_{\varepsilon}} }  ,
\]
\[
I_{8} = \gamma _{k} {\int_{C_{\varepsilon}}   {u_{k} \left( {x}
\right)x^{\left( {2\alpha}  \right)}\prod\limits_{i=1}^k\left(x_i\xi_i\right)^{1-2\alpha_i} \cdot r^{ - 2\beta _{k}} F_{A}^{\left(
{n} \right)} \left( {\beta _{k} ;\theta}  \right){\sum\limits_{i = 1}^{k}
{{\frac{{1 - 2\alpha _{i}}} {{x_{i}}} }\cos \left( {{\rm {\bf n}};x_{i}}
\right)}} dC_{\varepsilon}} }  .
\]

В интеграле $I_{5} $ переходим на обобщенную сферическую систему координат
$$x_{i} = \xi _{i} + \varepsilon \,\Phi _{i} \left( {\phi}  \right),i =
\overline {1,m} ,$$
где
\[
\phi: = \left( {\phi _{1} ,...,\phi _{n}}  \right),
\quad
\Phi _{1}: = \cos \phi _{1} ,
\Phi _{2}: = \sin \phi _{1} \cos \phi _{2} ,
\quad
\Phi _{3}: = \sin \phi _{1} \sin \phi _{2} \cos \phi _{3} ,
...,
\]
\[
\Phi _{m - 1}: = \sin \phi _{1} \sin \phi _{2} ...\sin \phi _{m - 2} \cos
\phi _{m - 1} ,
\quad
\Phi _{m}: = \sin \phi _{1} \sin \phi _{2} ...\sin \phi _{m - 2} \sin \phi
_{m - 1}
\]
\[
{\left[ {\varepsilon \ge 0,\,\,0 \le \phi _{1} \le \pi ,\,...,\,0 \le \phi
_{m - 2} \le \pi ,\,\,\,0 \le \phi _{m - 1} \le 2\pi \,\,} \right]}.
\]

Тогда мы имеем
\[
I_{5} = 2{\beta} _{k} \gamma _{k} \,\varepsilon ^{ - 2 {\beta}_{k} + m - 2}
\]
\[
\cdot{\int\limits_{0}^{2\pi}  {d\phi _{m - 1}}
}{\int\limits_{0}^{\pi}  {\sin \phi _{m - 2} d\phi _{m - 2}}
}...{\int\limits_{0}^{\pi}  {u_{k} \left( {\xi_{1} + \varepsilon \,\Phi
_{1} (\phi ),...,\xi _{m} + \varepsilon \,\Phi _{m} (\phi )} \right)}}
\]
\begin{equation*}
\label{eq65}
 \cdot {\prod\limits_{j = 1}^{k} {{\left[ {\xi _{j} + \varepsilon \,\Phi
_{j} (\phi )} \right]}}} {\prod\limits_{j = k + 1}^{n} {{\left[ {\left( {\xi
_{j} + \varepsilon \,\Phi _{j} (\phi )} \right)^{2\alpha _{j}}}  \right]}}
}F\left( {\varepsilon ,\phi}  \right)\sin ^{m - 2}\phi _{1} d\phi _{1} ,
\end{equation*}
где
\[
F\left( {\varepsilon ,\phi}  \right): = F_{A}^{(n)} {\left[
{{\begin{array}{*{20}c}
 {1 + \beta _{k} ,1 - \alpha _{1} ,...,1 - \alpha _{k} ,\alpha _{k + 1}
,...,\alpha _{n} ;} \hfill \\
 {2 - 2\alpha _{1} ,...,2 - 2\alpha _{k} ,2\alpha _{k + 1} ,...,2\alpha _{n}
;} \hfill \\
\end{array}} \Upsilon(\varepsilon ,\phi)} \right]},
\]
\[
\Upsilon(\varepsilon ,\phi):= \left(1 - {\frac{{\bar {r}_{1}^{2} \left( {\varepsilon ,\phi}
\right)}}{{\varepsilon ^{2}}}},...,1 - {\frac{{\bar {r}_{n}^{2} \left(
{\varepsilon ,\phi}  \right)}}{{\varepsilon ^{2}}}}\right),
\]
\[
\bar {r}_{i}^{2} \left( {\varepsilon ,\phi}  \right): = \left( {2\xi _{i} +
\varepsilon \,\Phi _{i} (\phi )} \right)^{2} + \varepsilon
^{2}{\sum\limits_{j = 1,j \ne i}^{m} {{\left[ {\Phi _{j} (\phi )\xi _{j}}
\right]}^{2}}} ,\,\,\,i = \overline {1,n} .
\]

Для полного вычисления $I_{3} $ сначала вычислим $F\left( {\varepsilon
,\phi}  \right)$. Используя последовательно формулу разложения (\ref{eq24}),
известную формулу
 \cite[c.113]{A2}:
\begin{equation}
\label{eq18}
F\left( {a,b;c;z} \right) = \left( {1 - z} \right)^{ - b}F\left( {c -
a,b;c;{\frac{{z}}{{z - 1}}}} \right),
\end{equation}
и формулу суммирования (\ref{eq17}),
получим
\[
 {\mathop {\lim} \limits_{\varepsilon \to 0}} \varepsilon ^{ - 2k + 2\alpha
_{1} + ... + 2\alpha _{k} - 2\alpha _{k + 1} - .. - 2\alpha _{n}} F\left(
{\varepsilon ,\phi}  \right) \\
\]
\[
 = {\prod\limits_{s = 1}^{k} {{\left[ {{\frac{{\Gamma \left( {2 - 2\alpha
_{s}}  \right)\Gamma \left( {\beta _{k} + \alpha _{s}}  \right)}}{{\Gamma
\left( {1 - \alpha _{s}}  \right)\Gamma \left( {\beta _{k} + 1}
\right)}}}\left( {2\xi _{s}}  \right)^{2\alpha _{s} - 2}} \right]}}
}{\prod\limits_{s = k + 1}^{n} {{\left[ {{\frac{{\Gamma \left( {2\alpha _{s}
} \right)\Gamma \left( {\beta _{k} - \alpha _{s} + 1} \right)}}{{\Gamma
\left( {\alpha _{s}}  \right)\Gamma \left( {\beta _{k} + 1} \right)}}}\left(
{2\xi _{s}}  \right)^{ - 2\alpha _{s}}}  \right]}}}  \\
\]
\[
{\sum\limits_{{\mathop {m_{i,j} =
0}\limits_{(2 \le i \le j \le n)}}} ^{\infty}  {{\frac{{\prod\limits_{s = 1}^{k} {\left[ {(1 - \alpha _{s}
)_{B(s,n)} \left( {\beta _{k} + \alpha _{s}}  \right)_{A(s,n) - B(s,n)}}
\right]} \prod\limits_{s = k + 1}^{n} {{\left[ {(\alpha _{s} )_{B(s,n)}
\left( {\beta _{k} - \alpha _{s} + 1} \right)_{A(s,n) - B(s,n)}}  \right]}}}}{{{\mathop {m_{2,2} ! \cdot \cdot \cdot m_{i,j} ! \cdot
\cdot \cdot m_{n,n} !}\prod\limits_{s = 1}^{n - 1} {{\left[ {\left( {\beta _{k} +
1} \right)_{A(s,n)}}  \right]}}}} }}}}.
\]

Далее,  с учетом формулы (\ref{eq28}), имеем
\[ {\mathop {\lim} \limits_{\varepsilon \to 0}} \varepsilon ^{ - 2k + 2\alpha
_{1} + ... + 2\alpha _{k} - 2\alpha _{k + 1} - .. - 2\alpha _{n}} F\left(
{\varepsilon ,\phi}  \right) \\
\]
\begin{equation}
\label{eq66}
 = {\prod\limits_{s = 1}^{k} {\left( {2\xi _{s}}  \right)^{2\alpha _{s} -
2}}} {\prod\limits_{s = k + 1}^{n} {\left( {2\xi _{s}}  \right)^{ - 2\alpha
_{s}}} } {\frac{{\Gamma \left( {m / 2} \right)}}{{\Gamma \left( {\beta _{k}
+ 1} \right)}}}{\prod\limits_{s = 1}^{k} {{\frac{{\Gamma \left( {2 - 2\alpha
_{s}}  \right)}}{{\Gamma \left( {1 - \alpha _{s}}  \right)}}}}
}{\prod\limits_{s = k + 1}^{n} {{\frac{{\Gamma \left( {2\alpha _{s}}
\right)}}{{\Gamma \left( {\alpha _{s}}  \right)}}}}} . \\
\end{equation}

Теперь переходим к пределу в $I_{5} $ при $\varepsilon \to 0$. В силу (\ref{eq32})
и (\ref{eq66})
окончательно находим
\begin{equation}
\label{eq67}
{\mathop {\lim} \limits_{\varepsilon \to 0}} I_{5} = u_{k} (\xi ).
\end{equation}

Аналогично следует, что
\begin{equation}
\label{eq68}
{\mathop {\lim} \limits_{\varepsilon \to 0}} I_{6} = {\mathop {\lim
}\limits_{\varepsilon \to 0}} I_{7} = {\mathop {\lim} \limits_{\varepsilon
\to 0}} I_{8} = {\mathop {\lim} \limits_{\varepsilon \to 0}} I_{4} =
{\mathop {\lim} \limits_{\varepsilon \to 0}} I_{2} = 0.
\end{equation}

Таким образом, левая часть равенства (\ref{eq58}) стала известной:
\begin{equation}
\label{eq69}
{\mathop {\lim} \limits_{\varepsilon \to 0}} I = {\mathop {\lim
}\limits_{\varepsilon \to 0}} {\int_{C_{\varepsilon}}   {x^{\left( {2\alpha
} \right)}{\left[ {u_{k} {\frac{{\partial G_{k}}} {{\partial {\rm {\bf
n}}}}} - G_{k} {\frac{{\partial u_{k}}} {{\partial {\rm {\bf n}}}}}}
\right]}dC_{\varepsilon}} }   = u_{k} \left( {\xi}  \right).
\end{equation}

Теперь займемся с правой частью равенства (\ref{eq58}). Рассмотрим интеграл
\begin{equation*}
\label{eq70}
{\int_{S} {\varphi _{k} \left( {S} \right)x^{\left( {2\alpha}
\right)}{\frac{{\partial G_{k} \left( {x;\xi}  \right)}}{{\partial {\rm {\bf
n}}}}}dS}} .
\end{equation*}

Поскольку на $S$ выполняются равенства $\bar {r}^{2} = r^{2},\bar
{r}_{1}^{2} = r_{1}^{2} $ и $\bar {\theta}  = \theta ,$ то, с учетом формул
(\ref{eq44}) и (\ref{eq55}), после нескольких элементарных преобразований, найдем
\begin{equation}
\label{eq71}
{\left. {{\frac{{\partial G_{k} \left( {x;{\rm \xi}}  \right)}}{{\partial
n}}}} \right|}_{S} = 2\beta _{k} \gamma _{k} \prod\limits_{i=1}^k\left(x_i\xi_i\right)^{1-2\alpha_i} \cdot F_{A}^{\left( {n} \right)}
\left( {1 + \beta _{k} ;\theta}  \right){\frac{{\rho ^{2} - R^{2}}}{{Rr^{2 +
2\beta _{k}}} }}.
\end{equation}

Подставив теперь (\ref{eq69}) и (\ref{eq71}) в формулу (\ref{eq58}), получим решение обобщенной
задачи Хольмгрена с условиями (\ref{eq13})-(\ref{eq15}) для уравнения (\ref{eq1}) в явном виде

\[
 u_{k} \left( {\xi}  \right) = {\sum\limits_{p = 1}^{k} {{\int_{D _{p}
} {\tau _{p} \left( {\tilde {x}_{p}}  \right)\tilde {G}_{k} \left(
{x_{p}^{0} ;\xi}  \right)d\tilde {x}_{p}}} } }  - {\sum\limits_{p = k +
1}^{n} {{\int_{D _{p}}  {\tilde {x}_{p}^{\left( {2\alpha}  \right)} \nu
_{p} \left( {\tilde {x}_{p}}  \right)G_{k} \left( {x_{p}^{0} ;\xi}
\right)d\tilde {x}_{p}}} } }  \\
\]
\begin{equation}
\label{eq72}
 + 2\beta _{k} \gamma _{k} {\int_{S} {\varphi _{k} \left( {x} \right)x^{\left(
{2\alpha}  \right)}\prod\limits_{i=1}^k\left(x_i\xi_i\right)^{1-2\alpha_i} \cdot F_{A}^{\left( {n} \right)} \left( {1 + \beta _{k} ;\theta}
\right){\frac{{R^{2} - \rho ^{2}}}{{Rr^{2 + 2\beta _{k}}} }}d_{x} S}} , \\
 \end{equation}
где
\[
\tilde {G}_{k} \left( {x_{i}^{0} ;\xi}  \right) = \left( {1 - 2\alpha _{i}}
\right)\gamma _{k} \tilde {x}_{i}^{(2\alpha )} \xi _{i}^{1 - 2\alpha _{i}}
\]
\begin{equation*}
\label{eq73}
\cdot{\prod\limits_{j = 1,j \ne i}^{k} {{\left[ {\left( {x_{j} \xi _{j}}
\right)^{1 - 2\alpha _{j}}}  \right]}}} \,{\left[ {{\frac{{F_{(i \le
k)}^{\left( {n - 1,0} \right)} \left( {\theta ^{\left( {0i} \right)}}
\right)}}{{r_{0i}^{2\beta _{k}}} } } - {\frac{{F_{(i \le k)}^{\left( {n -
1,0} \right)} \left( {\bar {\theta} ^{\left( {0i} \right)}} \right)}}{{\bar
{r}_{0i}^{2\beta _{k}}} } }} \right]},
\end{equation*}
\begin{equation*}
\label{eq74}
G_{k} \left( {x_{i}^{0} ;\xi}  \right) = \gamma _{k} \prod\limits_{i=1}^k\left(x_i\xi_i\right)^{1-2\alpha_i} \cdot {\left[ {{\frac{{F_{(i >
k)}^{\left( {n - 1,0} \right)} \left( {\theta ^{\left( {0i} \right)}}
\right)}}{{r_{0i}^{2\beta _{k}}} } } - {\frac{{F_{(i > k)}^{\left( {n -
1,0} \right)} \left( {\bar {\theta} ^{\left( {0i} \right)}} \right)}}{{\bar
{r}_{0i}^{2\beta _{k}}} } }} \right]},
\end{equation*}
\begin{equation*}
\label{eq751}
F_{(i \le k)}^{\left( {n - 1,0} \right)} \left( {z} \right)=F_{A}^{\left( {n - 1} \right)} {\left[ {{\begin{array}{*{20}c}
 \beta _{k} ,1 - \alpha _{1} ,..., \hfill \\
 2 - 2\alpha _{1} ,...,\hfill \\
\end{array}}} \right.}
\end{equation*}
\begin{equation*}
\label{eq75}
{\left. {{\begin{array}{*{20}c}
...,1 - \alpha _{i - 1}, 1 - \alpha _{i + 1}
,...,1 - \alpha _{k} ,\alpha _{k + 1} ,...,\alpha _{n} ; \hfill \\
...,2 - 2\alpha _{i - 1},2 - 2\alpha _{i + 1} ,...,2 -2\alpha _{k} ,2\alpha _{k + 1} ,...,2\alpha _{n} ; \hfill \\
\end{array}} z} \right]},
\end{equation*}
\begin{equation*}
\label{eq761}
F_{\left( {i > k} \right)}^{\left( {n - 1,0} \right)} \left( {z} \right): =
F_{A}^{\left( {n - 1} \right)} {\left[ {{\begin{array}{*{20}c}
 {\beta _{k} ,1 - \alpha _{1} ,...,} \hfill \\
 {2 - 2\alpha _{1} ,...,} \hfill \\
\end{array}} } \right.}
\end{equation*}
\begin{equation*}
\label{eq76}
{\left. {{\begin{array}{*{20}c}
 {...,1 - \alpha _{k} ,\alpha _{k + 1}
,...,\alpha _{i - 1} ,\alpha _{i + 1} ,...,\alpha _{n} ;} \hfill \\
 {...,2 - 2\alpha _{k} ,2\alpha _{k + 1} ,...,2\alpha _{i
- 1} ,2\alpha _{i + 1} ,...,2\alpha _{n} ;} \hfill \\
\end{array}} z} \right]},
\end{equation*}
\begin{equation*}
\label{eq77}
\theta ^{\left( {0i} \right)} = \left( {\theta _{1}^{\left( {0i} \right)}
,...,\theta _{i - 1}^{\left( {0i} \right)} ,\theta _{i + 1}^{\left( {0i}
\right)} ,...,\theta _{n}^{\left( {0i} \right)}}  \right),
\quad
\theta _{l}^{\left( {0i} \right)} = 1 - {\frac{{r_{0il}^{2}}} {{r_{0i}^{2}
}}},
\end{equation*}
\begin{equation*}
\label{eq771}
r_{0i}^{2} = \xi _{i}^{2} + {\sum\limits_{j = 1,j \ne i}^{m} {\left( {x_{j}
- \xi _{j}}  \right)^{2}}} ,
\end{equation*}
\begin{equation*}
\label{eq78}
\bar {r}_{0i}^{2} = {\sum\limits_{j = 1,j \ne i}^{m} {\left( {R -
{\frac{{x_{j} \xi _{j}}} {{R}}}} \right)^{2}}}  +
{\frac{{1}}{{R^{2}}}}{\sum\limits_{j = 1,j \ne i}^{m} {{\sum\limits_{l = 1,l
\ne j}^{m} {x_{j}^{2} \xi _{l}^{2}}}   - (m - 2)R^{2}}}
\end{equation*}
\begin{equation*}
\label{eq79}
r_{0il}^{2} = \xi _{i}^{2} + \left( {x_{l} + \xi _{l}}  \right)^{2} +
{\sum\limits_{j = 1,j \ne i,j \ne l}^{m} {\left( {x_{j} - \xi _{j}}
\right)^{2}}} ,
i,l = \overline {1,n} .
\end{equation*}

Отметим, что формула (\ref{eq72}) в некотором смысле обобщает известную формулу Пуассона задачи Дирихле для уравнения Лапласа (см., например, \cite[c. 269]{Mix}), поэтому формула (\ref{eq72})  называется \textit{сингулярной формулой Пуассона} (или \textit{сингулярным аналогом формулы Пуассона}), а выражение
\[
P\left( {x;\xi}  \right) = x^{\left( {2\alpha}  \right)}\prod\limits_{i=1}^k\left(x_i\xi_i\right)^{1-2\alpha_i} \cdot
\]
\begin{equation}
\label{eq80}
 \cdot F_{A}^{\left( {n} \right)}
{\left[ {{\begin{array}{*{20}c}
 {1 + \beta _{k} ,1 - \alpha _{1} ,...,1 - \alpha _{k} ,\alpha _{k + 1}
,...,\alpha _{n} ;} \hfill \\
 {2 - 2\alpha _{1} ,...,2 - 2\alpha _{k} ,2\alpha _{k + 1} ,...,2\alpha _{n}
;} \hfill \\
\end{array}} \theta}  \right]}{\frac{{R^{2} - \rho ^{2}}}{{Rr^{2 + 2\beta
_{k}}} }},
\quad
\rho \le R,
\end{equation}
-- \textit{сингулярным ядром Пуассона} (или \textit{сингулярным аналогом ядра Пуассона}).

Покажем, что функция определенная формулой (\ref{eq72}) действительно является
решением обобщенной задачи Холмгрена. С этой целью решение (\ref{eq72}) перепишем в
виде
\begin{equation}
\label{eq81}
u_{k} \left( {\xi}  \right) = {\sum\limits_{l = 1}^{k} {{\left[ {I_{l}
\left( {\xi}  \right) + J_{l} \left( {\xi}  \right)} \right]}}}  +
{\sum\limits_{l = k + 1}^{n} {{\left[ {K_{l} \left( {\xi}  \right) + L_{l}
\left( {\xi}  \right)} \right]}}}  + M\left( {\xi}  \right),
\end{equation}
где
\[
I_{l} \left( {\xi}  \right) = \left( {1 - 2\alpha _{l}}  \right)\gamma _{k} \prod\limits_{j=1}^k\left[\xi_j^{1-2\alpha_j}\right]\int_{D _{l}}  \tau _{l}
\left( {\tilde {x}_{l}}  \right)
\]
\begin{equation}
\label{eq82}
 \cdot  {{\prod\limits_{j
= 1,j \ne l}^{k} {{\left[ {x_{j}}   \right]}}}{\prod\limits_{j
= k+1}^{n} {{\left[ {x_{j}^{2\alpha _{j}}}   \right]}}} \,r_{0l}^{
- 2\beta _{k}}  F_{(l \le k)}^{\left( {n - 1,0} \right)} \left( {\theta
^{\left( {0l} \right)}} \right)\,d\tilde {x}_{l}}  ,\,\,l = \overline {1,k}
,
\end{equation}
\[
J_{l} \left( {\xi}  \right) = - \left( {1 - 2\alpha _{l}}  \right)\gamma
_{k} \prod\limits_{j=1}^k\left[\xi_j^{1-2\alpha_j}\right]\int_{D_{l}}  \tau
_{l} \left( {\tilde {x}_{l}}  \right)
\]
\[
\cdot {{\prod\limits_{j
= 1,j \ne l}^{k} {{\left[ {x_{j}}   \right]}}}{\prod\limits_{j
= k+1}^{n} {{\left[ {x_{j}^{2\alpha _{j}}}   \right]}}} \,\bar {r}_{0l}^{ - 2\beta _{k}}  F_{(l \le k)}^{\left( {n -
1,0} \right)} \left( {\bar {\theta} ^{\left( {0l} \right)}} \right)\,d\tilde
{x}_{l}}  ,\,\,l = \overline {1,k} ,
\]
\[
K_{l} \left( {\xi}  \right) = - \gamma _{k} \prod\limits_{j=1}^k\left[\xi_j^{1-2\alpha_j}\right]\int_{D _{l}}  \nu _{l} \left( {\tilde {x}_{l}}
\right)
\]
\[
\cdot {{\prod\limits_{j=1,j \ne l}^{k} {{\left[ {x_{j}}   \right]}}}{\prod\limits_{j
= k+1}^{n} {{\left[ {x_{j}^{2\alpha _{j}}}   \right]}}} r_{0l}^{ - 2\beta _{k}}  F_{(l > k)}^{\left( {n - 1,0} \right)}
\left( {\theta ^{\left( {0l} \right)}} \right)d\tilde {x}_{l}}  ,\,\,l =
\overline {k + 1,n} ,
\]
\[
L_{l} \left( {\xi}  \right) = \gamma _{k} \prod\limits_{j=1}^k\left[\xi_j^{1-2\alpha_j}\right]\int_{D_{l}}  \nu _{l} \left( {\tilde {x}_{l}}
\right)
\]
\[
 \cdot {{\prod\limits_{j= 1,j \ne l}^{k} {{\left[ {x_{j}}   \right]}}}{\prod\limits_{j
= k+1}^{n} {{\left[ {x_{j}^{2\alpha _{j}}}   \right]}}} \bar {r}_{0l}^{ - 2\beta _{k}}  F_{(l > k)}^{\left( {n - 1,0}
\right)} \left( {\bar {\theta} ^{\left( {0l} \right)}} \right)d\tilde
{x}_{l}}  ,\,\,\,l = \overline {k + 1,n} ,
\]
\[
M\left( {\xi}  \right) = 2\beta _{k} \gamma _{k} \prod\limits_{j=1}^k\left[\xi_j^{1-2\alpha_j}\right] \int_{S} \varphi _{k} \left( {S} \right)
\]
\[
\cdot {\prod\limits_{j
= 1}^{k} {{\left[ {x_{j}}   \right]}}}{\prod\limits_{j
= k+1}^{n} {{\left[ {x_{j}^{2\alpha _{j}}}   \right]}}} F_{A}^{\left( {n} \right)}
\left( {1 + \beta _{k} ;\theta}  \right){\frac{{R^{2} - \rho ^{2}}}{{Rr^{2 +
2\beta _{k}}} }}d_{x} S,
\]

Вычислим предел ${\mathop {\lim} \limits_{\xi _{s} \to 0}} I_{l} \left( {\xi
} \right)$ как при $s = l,$ так и при $s \ne l\left( {l,s = \overline
{1,k}}  \right).$

Предварительно преобразуем функцию $F_{(l \le k)}^{\left( {n - 1,0}
\right)} \left( {\theta ^{\left( {0l} \right)}} \right)$, входящую в
интеграл $I_{l} \left( {\xi}  \right)$. В силу следствия 1, получим
\[
 F_{(l \le k)}^{\left( {n - 1,0} \right)} \left( {\theta ^{\left( {0l}
\right)}} \right) ={\sum\limits_{{\mathop {m_{i,j} =
0}\limits_{(2 \le i \le j \le n-1)}}} ^{\infty}  {{\frac{{(\beta_k)_{A(n-1)}}}{{{\mathop {m_{2,2} !m_{2,3} ! \cdot \cdot \cdot m_{i,j} ! \cdot
\cdot \cdot m_{n-1,n-1} !}}} }}}}
\\
\]
\[
 \cdot {\prod\limits_{s = 1}^{l - 1} {{ {{\frac{{(1 - \alpha _{s}
)_{B(s)}}} {{(2 - 2\alpha _{s} )_{B(s)}}} }\left( {\theta _{s}^{\left( {0l}
\right)}}  \right)^{B(s)}F\left( {{\begin{array}{*{20}c}
 {\beta _{k} + A(s),1 - \alpha _{s} + B(s);} \hfill \\
 {2 - 2\alpha _{s} + B(s);} \hfill \\
\end{array}} \theta _{s}^{\left( {0l} \right)}}  \right)}}}}  \\
\]
\[
 \cdot {\prod\limits_{s = l + 1}^{k} {{{{\frac{{(1 - \alpha _{s}
)_{B(s - 1)}}} {{(2 - 2\alpha _{s} )_{B(s - 1)}}} }\left( {\theta
_{s}^{\left( {0l} \right)}}  \right)^{B(s - 1)}F\left(
{{\begin{array}{*{20}c}
 {\beta _{k} + A(s - 1),1 - \alpha _{s} + B(s - 1);} \hfill \\
 {2 - 2\alpha _{s} + B(s - 1);} \hfill \\
\end{array}} \theta _{s}^{\left( {0l} \right)}}  \right)}}}}  \\
\]
\begin{equation}
\label{eq83}
 \cdot {\prod\limits_{s = k + 1}^{n} {{{{\frac{{(\alpha _{s} )_{B(s -
1)}}} {{(2\alpha _{s} )_{B(s - 1)}}} }\left( {\theta _{s}^{\left( {0l}
\right)}}  \right)^{B(s - 1)}F\left( {{\begin{array}{*{20}c}
 {\beta _{k} + A(s - 1),\alpha _{s} + B(s - 1);} \hfill \\
 {2\alpha _{s} + B(s - 1);} \hfill \\
\end{array}} \theta _{s}^{\left( {0l} \right)}}  \right)}}}} . \\
\end{equation}

Здесь и далее, для краткости принята сокращенная запись:
\[
A(n - 1): = A(n - 1,n - 1),\,
A(s): = A(s,n - 1),\,
B(s): = B(s,n - 1);
\]
\[
A(s - 1): = A(s - 1,n - 1),\,
B(s): = B(s - 1,n - 1).
\]

Применение формулы (\ref{eq18}) к каждой гипергеометрической функции Гаусса,
входящей в формулу (\ref{eq83}), дает
\begin{equation}
\label{eq84}
F_{(l \le k)}^{\left( {n - 1,0} \right)} \left( {\theta ^{\left( {0l}
\right)}} \right) = {\prod\limits_{s = 1,s \ne l}^{k} {\left(
{{\frac{{r_{0l}^{2}}} {{r_{0ls}^{2}}} }} \right)^{1 - \alpha _{s}}}
}{\prod\limits_{s = k + 1}^{n} {\left( {{\frac{{r_{0l}^{2}}} {{r_{0ls}^{2}
}}}} \right)^{\alpha _{s}}} } \tilde {F}_{(l \le k)}^{\left( {n - 1,0}
\right)} \left( Z \right),
\end{equation}
где
\[
\tilde {F}_{(l \le k)}^{\left( {n - 1,0} \right)} \left(Z \right) = {\sum\limits_{{\mathop {m_{i,j} =
0}\limits_{(2 \le i \le j \le n-1)}}} ^{\infty}  {{\displaystyle\frac{{(\beta_k)_{A(n-1)}}}{{{\mathop {m_{2,2}!\cdot \cdot \cdot m_{i,j}!\cdot \cdot \cdot m_{n-1,n-1}!}}} }}}}{\prod\limits_{s =
1}^{l - 1}{\displaystyle\frac{{(1 - \alpha _{s}
)_{B(s)}}} {{(2 - 2\alpha _{s} )_{B(s)}}} } {\left(-Z\right)^{B(s)}}}
\]
\[
\cdot {\prod\limits_{s = l + 1}^{k}{\displaystyle\frac{{(1 - \alpha _{s}
)_{B(s - 1)}}} {{(2 - 2\alpha _{s} )_{B(s - 1)}}} } {\left(-Z \right)^{B(s - 1)}}} {\prod\limits_{s = k + 1}^{n}{\displaystyle\frac{{(\alpha _{s} )_{B(s -
1)}}} {{(2\alpha _{s} )_{B(s - 1)}}} }
{\left( -Z \right)^{B(s - 1)}}}
\]
\[
 \cdot {\prod\limits_{s = 1}^{l - 1} {{ {F\left( {{\begin{array}{*{20}c}
 {2 - 2\alpha _{s} - \beta _{k} + B(s) - A(s),1 - \alpha _{s} + B(s);}
\hfill \\
 {2 - 2\alpha _{s} + B(s);} \hfill \\
\end{array}} Z} \right)} }}
} \]
\[
 \cdot {\prod\limits_{s = l + 1}^{k} {{ {F\left(
{{\begin{array}{*{20}c}
 {2 - 2\alpha _{s} - \beta _{k} + B(s - 1) - A(s - 1),1 - \alpha _{s} + B(s
- 1);} \hfill \\
 {2 - 2\alpha _{s} + B(s - 1);} \hfill \\
\end{array}} Z} \right)} }}
} \]
\[
 \cdot {\prod\limits_{s = k + 1}^{n} {{ {F\left(
{{\begin{array}{*{20}c}
 {2\alpha _{s} - \beta _{k} + B(s - 1) - A(s - 1),\alpha _{s} + B(s - 1);}
\hfill \\
 {2\alpha _{s} + B(s - 1);} \hfill \\
\end{array}} Z} \right)} }}
},\]
\[\,\,Z:={1-{\displaystyle\frac{{r_{0l}^{2}}} {{r_{0ls}^{2}}}}}.\]

Вычислим значение функции $\tilde {F}_{(l \le k)}^{\left( {n - 1,0}
\right)} \left( Z \right)$ при $r_{0l}^{2} \to 0$.
После применения формулы суммирования (\ref{eq17}), получим
\begin{equation}
\label{eq86}
\begin{array}{l}
 {\mathop {\lim} \limits_{r_{0l} \to 0}} \tilde {F}_{(l \le k)}^{\left( {n
- 1,0} \right)} \left(Z\right) = {\prod\limits_{s
= 1,s \ne l}^{k} {{ {{\displaystyle\frac{{\Gamma \left( {\beta _{k} - 1 +
\alpha _{s}}  \right)\Gamma \left( {2 - 2\alpha _{s}}  \right)}}{{\Gamma
\left( {\beta _{k}}  \right)\Gamma \left( {1 - \alpha _{s}}  \right)}}}}
}}} {\prod\limits_{s = k + 1}^{n} {{ {{\displaystyle\frac{{\Gamma \left(
{\beta _{k} - \alpha _{s}}  \right)\Gamma \left( {2\alpha _{s}}
\right)}}{{\Gamma \left( {\beta _{k}}  \right)\Gamma \left( {\alpha _{s}}
\right)}}}}}}}  \\
\\
 \cdot {\sum\limits_{{\mathop {m_{i,j} =
0}\limits_{(2 \le i \le j \le n-1)}}} ^{\infty}  {{\displaystyle\frac{{(\beta_k)_{A(n-1)}}}{{{\mathop {m_{2,2} ! \cdot \cdot \cdot m_{i,j} ! \cdot\cdot \cdot m_{n-1,n-1} !}}} }}}}{\prod\limits_{s = 1}^{l - 1} {{{{\displaystyle\frac{{(1 -\alpha _{s} )_{B(s)} \left( {\beta _{k} - 1 + \alpha _{s}}  \right)_{A(s) -
B(s)}}} {{\left( {\beta _{k}}  \right)_{A(s)}}} }}}}}  \\
\\
 \cdot {\prod\limits_{s = l + 1}^{k} {{ {{\displaystyle\frac{{(1 - \alpha _{s}
)_{B(s - 1)} \left( {\beta _{k} - 1 + \alpha _{s}}  \right)_{A(s - 1) - B(s
- 1)}}} {{\left( {\beta _{k}}  \right)_{A(s - 1)}}} }}}}
}{\prod\limits_{s = k + 1}^{n} {{ {{\displaystyle\frac{{(\alpha _{s} )_{B(s - 1)}
\left( {\beta _{k} - \alpha _{s}}  \right)_{A(s - 1) - B(s - 1)}}} {{\left(
{\beta _{k}}  \right)_{A(s - 1)}}} }}}}} . \\
 \end{array}
\end{equation}

Пользуясь теперь формулой (\ref{eq27}), окончательно получим
\begin{equation}
\label{eq87}
{\mathop {\lim} \limits_{r_{0l} \to 0}} \tilde {F}_{(l \le k)}^{\left( {n -
1,0} \right)} \left( Z  \right) = \Gamma ^{ -
1}\left( {\beta _{k}}  \right)\Gamma \left( {{\frac{{m}}{{2}}} - \alpha _{l}
} \right){\prod\limits_{s = 1,s \ne l}^{k} {{ {{\frac{{\Gamma \left(
{2 - 2\alpha _{s}}  \right)}}{{\Gamma \left( {1 - \alpha _{s}}  \right)}}}}
}}} {\prod\limits_{s = k + 1}^{n} {{ {{\frac{{\Gamma \left(
{2\alpha _{s}}  \right)}}{{\Gamma \left( {\alpha _{s}}  \right)}}}}}}} .
\end{equation}

Теперь вычислим предел ${\mathop {\lim} \limits_{\xi _{s} \to 0}} I_{l}
\left( {\xi}  \right)$ при $s = l.$ Действительно, положив в правой
части равенства (\ref{eq82})
\begin{equation}
\label{eq88}
x_{i} = \xi _{i} + \xi _{l} t_{i} ,\,\,\,i = \overline {1,m} ,\,\,i \ne l
\end{equation}
переходим в $I_{s} \left( {\xi}  \right)$ к пределу при $\xi _{s} \to 0$ .
Применяя равенство (\ref{eq87}) и учитывая выражение (\ref{eq32}) для $\gamma _{k} $,
получим
\[
{\mathop {\lim} \limits_{\xi _{l} \to 0}} I_{l} \left( {\xi}  \right) = 2^{
- 2\alpha _{l}} \pi ^{ - m / 2}\Gamma \left( {{\frac{{m}}{{2}}} - \alpha
_{l}}  \right)
\]
\begin{equation}
\label{eq89}
\cdot {\frac{{\Gamma \left( {1 - \alpha _{l}}  \right)}}{{\Gamma
\left( {1 - 2\alpha _{l}}  \right)}}}{\int\limits_{ - \infty} ^{\infty}
{...{\int\limits_{ - \infty} ^{\infty}  {\left( {1 + {\sum\limits_{i = 1,i
\ne l}^{m} {t_{i}^{2}}} }  \right)^{\alpha _{l} - {\frac{{m}}{{2}}}}d\tilde
{t}_{l}}} } }  \cdot \tau \left( {\tilde {x}_{l}}  \right),
\end{equation}
где
$$d\tilde t_l:= dt_1\cdot\cdot\cdot dt_{l-1}dt_{l+1}\cdot\cdot\cdot dt_m.$$

Учитывая \cite[c.637]{A24}
\[
{\int\limits_{ - \infty} ^{\infty}  {...{\int\limits_{ - \infty} ^{\infty}
{\left( {1 + {\sum\limits_{i = 1,i \ne l}^{m} {t_{i}^{2}}} }
\right)^{\displaystyle\alpha _{l} - {\displaystyle\frac{{m}}{{2}}}}d\tilde {t}_{l}}} } }  =
{\displaystyle\frac{{\Gamma ^{m - 1}\left( {{\displaystyle\frac{{1}}{{2}}}} \right)\Gamma \left(
{{\displaystyle\frac{{1}}{{2}}} - \alpha _{l}}  \right)}}{{\Gamma \left(
{{\displaystyle\frac{{m}}{{2}}} - \alpha _{l}}  \right)}}},
\]
из (\ref{eq89}) получим, что
\begin{equation}
\label{eq90}
{\mathop {\lim} \limits_{\xi _{l} \to 0}} I_{l} \left( {\xi}  \right) = \tau
_{l} \left( {\tilde {x}_{l}}  \right),\,\,\,l = \overline {1,k} .
\end{equation}

Аналогичными вычислениями можно показать, что
\begin{equation}
\label{eq91}
{\mathop {\lim} \limits_{\xi _{s} \to 0}} I_{l} \left( {\xi}  \right) =
0,\,\,\,l,s = \overline {1,k} ,\,\,l \ne s;
\quad
{\mathop {\lim} \limits_{\xi _{s} \to 0}} J_{l} \left( {\xi}  \right) =
0,\,\,\,l,s = \overline {1,k} ;
\end{equation}
\begin{equation}
\label{eq92}
{\mathop {\lim} \limits_{\xi _{s} \to 0}} K_{l} \left( {\xi}  \right) =
{\mathop {\lim} \limits_{\xi _{s} \to 0}} L_{l} \left( {\xi}  \right) =
0,\,\,s = \overline {1,k} ,\,\,\,l = \overline {k + 1,n} ;
\quad
{\mathop {\lim} \limits_{\xi _{s} \to 0}} M\left( {\xi}  \right) = 0,\,\,\,s
= \overline {1,k} .
\end{equation}

Таким образом, функция, определенная формулой (\ref{eq72}), удовлетворяет условиям
(\ref{eq13}). Аналогично, можно показать выполнение условий (\ref{eq14}) и (\ref{eq15}).

Формула (\ref{eq72}), а с ней и все доказательство, требует, чтобы $m > 2$ . Однако сингулярная
формула Пуассона верна и для $m = 2$.

Таким образом, доказана

\textbf{Теорема 2}. \textit{Сингулярная формула Пуассона (\ref{eq72}) в области $\Omega$ является единственным решением обобщенной задачи Хольмгрена для уравнения (\ref{eq1}) с условиями (\ref{eq13})--(\ref{eq15})}.

\end{document}